\documentclass[10pt,onesite,draft]{article}
\usepackage{amsmath,amsxtra,amssymb,latexsym, amsfonts, indentfirst}
\usepackage[mathscr]{eucal}
\newfam\cyrfam

\newenvironment{proof}{%
\par\addvspace{6pt plus3pt minus2pt}%
\noindent{\bfseries\itshape\textit Proof.}\ignorespaces} {%
\if@halmos\halmos\fi
\par\addvspace{6pt plus3pt minus2pt} }
\begin{document}
\parskip 4pt
\large
\setlength{\baselineskip}{15 truept}
\setlength{\oddsidemargin} {0.5in}
\overfullrule=0mm
\def\bfh{\vhtimeb}
\date{}

\title{\bf \large A SURVEY OF SOME RESULTS \\ FOR MIXED
MULTIPLICITIES\\}
                                                                                                                                                                                                                                    \def\b{\vntime}
\author{
 Le Van Dinh  and  Nguyen Tien Manh\\ Truong Thi Hong Thanh \\
\small Department of Mathematics, Hanoi National University of Education\\
\small 136 Xuan Thuy Street, Hanoi, Vietnam\\
\small dinhlevands@gmail.com, manhnt79@gmail.com, thanhtth@hnue.edu.vn\\}
  \date{}
\maketitle
\centerline{
\parbox[c]{11. cm}{
\small  ABSTRACT: The paper proves that the superficial sequences
in \cite{Te, Tr2, TV,  VM} are the weak-(FC)-sequences in
\cite{Vi1}. This follows that
 \cite[Theorem 3.4]{Vi1} of Viet in 2000 covers the results of Risler and Teissier \cite{Te} in 1973;
  Trung \cite[Theorem 3.4]{Tr2} in 2001;          Trung and Verma \cite[Theorem 1.4]{TV} in 2007. The paper unified the results for mixed multiplicities of ideals
in \cite{Te, Tr2, TV,TV2,Vi1, VM}.}} \vspace{12pt}
\centerline{\Large\bf 1. Introduction}$\\$ In past years, using
different sequences, one expressed
  mixed multiplicities in terms of the Hilbert-Samuel multiplicity, for instance,
  in the case of ideals of
dimension $0$, Risler and Teissier \cite{Te} in 1973 showed that
each mixed multiplicity is the multiplicity of an ideal generated
by a superficial sequence. However, for the case of arbitrary
ideals, how to  find mixed multiplicity formulas, which are
analogous to Risler-Teissier's formula in \cite{Te}, was a
challenge. And this problem became an open question of the mixed
multiplicity theory (see e.g. \cite{KV, Tr0}). Viet \cite{Vi1} in
2000 solved this open problem \cite[Theorem 3.4]{Vi1}  via
(FC)-sequences. \footnotetext{\begin{itemize} \item[ ]{\it
Mathematics  Subject Classification} (2010): Primary 13H15.
Secondary 13C13; 13C15; 14C17; 13D40. \item[ ]$ Key\; words \; and
\; phrases:$ Noetherian ring, mixed multiplicity, superficial
sequence; weak-(FC)-sequence.
\end{itemize}}
Trung \cite{Tr2} in 2001 and Trung  and  Verma \cite{TV} in 2007
obtained  results that are similar   to the result  of \cite{Vi1}
with the help of  ``bi-filter-regular sequences'' and
$(\varepsilon_1,\ldots,\varepsilon_m)$-superficial sequences,
respectively.
\enlargethispage{0.5cm}
This paper
 seems to
account well for the minimum of conditions of weak-(FC)-sequences that is used as a tool in transmuting  mixed multiplicities of ideals into the Hilbert-Samuel multiplicity. This explains why the superficial sequences in \cite{Te, Tr2, TV, VM} are weak-(FC)-sequences; and
 \cite[Theorem 3.4]{Vi1} of Viet in 2000 covers the results of Risler   and  Teissier \cite{Te} in 1973; Trung \cite[Theorem 3.4]{Tr2} in 2001;   Trung  and  Verma \cite[Theorem 1.4]{TV} in 2007.

By studying \cite{Tr2, TV, Vi1},  we see that one can use Viet's
results and method in \cite{Vi1} for \cite{Tr2, TV} with replacing
weak-(FC)-sequences by other sequences that are variant forms of
weak-(FC)-sequences. In fact, this paper shows that the
superficial sequences in \cite{Te, Tr2, TV,  VM} are the
weak-(FC)-sequences in \cite{Vi1}. This proves that \cite[Theorem
3.4]{Tr2} of Trung in 2001 and \cite[Theorem 1.4]{TV} of Trung
and  Verma in 2007 are only variants  of \cite[Theorem 3.4]{Vi1}
of Viet in 2000. So the paper unified the results for mixed
multiplicities of ideals in \cite{Te, Tr2, TV,TV2,Vi1, VM}.

\vspace*{16pt}

\centerline{\Large\bf2.  Mixed Multiplicities  } \vspace*{8pt}

\noindent In this section, we introduce the notion of mixed
multiplicities of ideals and recall some results for mixed
multiplicities in \cite{Vi1}.

Let $(A, \frak m)$  be  a  noetherian   local ring with  maximal ideal $\frak{m}$ and the infinite residue field $k = A/\frak{m}$. Let $M$ be  a finitely generated  $A$-module. Let $J, I_1,\ldots, I_s$ be ideals of $A$  such that   $J$ is  an $\frak m$-primary ideal and $I=I_1\cdots I_s$  is not contained in $ \sqrt{\mathrm{Ann}{M}}$.
Set $\dim M/0_M: I^{\infty} = q.$

Then Viet  \cite [Proposition 3.1(i)]{Vi1} in 2000 (see \cite[Proposition 3.1(i)]{MV})
 proved that
$$\ell_A\bigg[\dfrac{J^{n_0}I_1^{n_1}\cdots I_s^{n_s}{M}}{J^{n_0+1}I_1^{n_1}\cdots I_s^{n_s}{M}}\bigg]$$
is a polynomial of degree $q-1$ for all large $n_0, n_1,\ldots,n_s.$ The terms of total degree $q-1$ in this polynomial have the form
$$ \sum_{k_0\:+k_1\:\cdots\:+\:k_s\;=\;q-1}e_A(J^{[k_0+1]},I_1^{[k_1]},\ldots,I_s^{[k_s]}; M)\dfrac{n_0^{k_0}n_1^{k_1}\cdots n_s^{k_s}}{k_0!k_1!\cdots k_s!}.$$
Then $e_A(J^{[k_0+1]},I_1^{[k_1]},\ldots,I_s^{[k_s]}; M)$ is
called the {\it  mixed multiplicity of $M$ with respect to ideals
$J,I_1,\ldots,I_s$ of the type $(k_0, k_1,\ldots,k_s).$} And Viet
\cite [Lemma 3.2(i)]{Vi1} in 2000 (see \cite[Lemma 3.2(i)]{MV})
showed that
$$e_A(J^{[k_0+1]},I_1^{[0]},\ldots,I_s^{[0]}; M) = e_A(J; M/0_M:
I^{\infty}).$$ \vskip 0.2cm \enlargethispage{3cm} \noindent {\bf
Remark 2.1.}  Trung  and  Verma \cite{TV2} in 2010 called (see
\cite [Lemma 5.1]{TV2}) that \cite [Proposition 3.1(i) and Lemma
3.2(i)]{Vi1} are
  results of Trung \cite{Tr2} in 2001 in the case when $s=1,$ and
are also  results of Trung  and  Verma in 2007 (see \cite [Theorem
8.1]{TV2}). Moreover, Trung  and  Verma \cite{TV} in 2007
considered \cite [Proposition 3.1(i) and  Lemma 3.2(i)]{Vi1}  as
Trung-Verma's results (see \cite [Theorem 1.2]{TV}). Although
\cite[Lemma 3.1]{Tr2} of Trung in 2001 is a particular case of
\cite [Proposition 3.1(i) and Lemma 3.2(i)]{Vi1} in 2000 but the
author of \cite{Tr2} seemed to omit these results.

\newpage
\vspace*{24pt} \centerline{\Large\bf3. Superficial Sequences
 } \vspace*{8pt}

\noindent
This section  proves  that the mentioned sequences in \cite{Te, Tr2, TV, Vi1,  VM} are weak-(FC)-sequences. As an application, we unify the results that interpreted mixed multiplicities of ideals as the Hilbert-Samuel multiplicity in the above works.

 Now, we recall superficial sequences which were chosen in Risler-Teissier's work \cite{Te}
 in 1973 (see e.g. \cite{CP, HS, K, RV, ZS}).

\vskip 0.2cm
\noindent
{\bf Definition 3.1.} An element $x \in I_i$
 is called an {\it $I_i$-superficial element of $M$ with respect to }$( I_1,\ldots,I_s)$ if  there exists a non-negative integer $c$ such that
$$(I_1^{n_1}\cdots I_i^{n_i+1}\cdots I_s^{n_s}M: x)\bigcap
I_1^{n_1}\cdots I_i^{c}\cdots I_s^{n_s}M = I_1^{n_1}\cdots I_s^{n_s}M$$
for all $n_i \ge c$ and for all non-negative integers $ n_1, \ldots, n_{i - 1}, n_{i + 1}, \ldots, n_s.$

Next, we recall the following sequence in \cite{Vi1} (see e.g.
\cite{CP1, DV, DV1, MV, ViFC, Vi2, Vi3,  VT}). \vskip 0.2cm
  \noindent {\bf Definition 3.2} \cite{Vi1}. Let $M$ be a finitely generated  $A$-module. Let $I_1,\ldots,I_s$ be ideals such that $I=I_1\cdots I_s$  is not contained in $ \sqrt{\mathrm{Ann}{M}}$.    An element $x \in A$ is called a {\it weak-$(FC)$- element of $M$ with respect to $(I_1,\ldots, I_s)$}  if there exists $i \in \{ 1, \ldots, s\}$ such that $x \in I_i$ and
the following conditions are satisfied:
 \begin{itemize}

      \item[$\mathrm{(i)}$] $x$ is an $I$-filter-regular element with respect to $M,$ i.e.,\;$0_M:x \subseteq 0_M: I^{\infty}.$
    \item[$\mathrm{(ii)}$] $xM\bigcap {I_1}^{n_1} \cdots I_i^{n_i+1}\cdots I_s^{n_s}M
= x{I_1}^{n_1}\cdots I_i^{n_i}\cdots I_s^{n_s}M$
for all $n_1,\ldots,n_s\gg0.$
 \end{itemize}
\noindent
Let $x_1, \ldots, x_t$ be a sequence in $A$. For any $0\le i < t,$ set\;      ${M}_i = \dfrac{M}{(x_1, \ldots, x_{i})M}$.  Then
$x_1, \ldots, x_t$ is called a {\it weak-$(FC)$-sequence of $M$ with respect to $(I_1,\ldots, I_s)$} if $x_{i + 1}$ is a weak-(FC)-element of ${M}_i$ with respect to $(I_1,\ldots, I_s)$  for all $i = 0, \ldots, t - 1$.
\vskip 0.2cm

Using weak-(FC)-sequences, the author of \cite{Vi1} in 2000
characterized mixed multiplicities as the Hilbert-Samuel
multiplicity (see e.g. \cite{DV, DV1, MV,Vi2, Vi3,  VT}) by the
following theorem. \enlargethispage{1cm}

\noindent{\bf Theorem 3.3} \cite[Theorem 3.4]{Vi1}. {\it Let $M$
be a finitely generated  $A$-module. Let $J, I_1,\ldots,I_s$ be
ideals of $A$ with $J$ being $\frak m$-primary and $I=I_1\cdots
I_s \nsubseteq \sqrt{\mathrm{Ann}{M}}$. Then
$e_A(J^{[k_0+1]},I_1^{[k_1]},\ldots,I_s^{[k_s]}; M) \ne 0$ if and
only if there exists a weak-$(FC)$-sequence $x_1,\ldots,x_m $ of
$M$ with respect to $(J, I_1,\ldots, I_s)$ consisting of $k_1$
elements of $I_1,\ldots,$ $k_s$ elements of $I_s$ such that $\dim
M/(x_1,\ldots,x_m)M:I^\infty = \dim M/0_M:I^\infty-m.$ In this
case,
$$e_A(J^{[k_0+1]},I_1^{[k_1]},\ldots,I_s^{[k_s]}; M)= e\big(J, \dfrac{M}{(x_1,\ldots,x_m)M:I^\infty}\big).$$}

For any graded algebra $S=\bigoplus_{n_1,\ldots,n_d\ge 0}
S_{(n_1,\ldots,n_d)},$ set  $ S_{++}=\bigoplus_{\;n_1,\ldots,n_d>
0}S_{(n_1,\ldots,n_d)}.$

\vskip 0.2cm  {\bf Note 3.4:}
 We recall sequences in \cite{Tr2,
TV} as follows.
\begin{itemize}
\item[$\mathrm{(i)}$]
In the case of two ideals, Trung in 2001 gave a result \cite[Theorem 3.4]{Tr2} similar   to Theorem 3.3 via ``bi-filter-regular sequences''. Set
$$F_J(J,I_1; A) = \bigoplus_{m,n\ge 0}\frac{J^mI_1^n}{J^{m+1}I_1^n}\;\; \mathrm{ and}\;\; F_{I_1}(I_1, J; A) = \bigoplus_{m,n\ge 0}\frac{J^mI_1^n}{J^{m}I_1^{n+1}}.$$    For any $x \in I_1,$ denote by $x^*$ and $x^{**}$ the images of $x$ in $I_1/JI_1$ and $I_1/I_1^2,$ respectively.  Let  $x_1,\ldots,x_t$ be a sequence in $I_1$. For any $1\le i < t,$\; set      ${A}_i = \dfrac{A}{(x_1, \ldots, x_{i})};$  $J_i = JA_i; {I_1}_i = I_1A_i$ and $\bar{x}_{i+1}$ is the image of $x_{i+1}$ in ${I_1}_i.$
Then Trung in \cite{Tr2} used a sequence $x_1,\ldots,x_t$ in $I_1$ such that
$x_1^*$ is an $F_{J}(J, I_1; A)_{++}$-filter-regular element and $x_1^{**}$ is an $
F_{I_1}(I_1,J; A)_{++}$-filter-regular element; $\bar{x}_{i+1}^*$ is an $F_{J_i}(J_i, {I_1}_i; A)_{++}$-filter-regular element and
$\bar{x}_{i+1}^{**}$ is an $F_{{I_1}_i}({I_1}_i,J_i; A)_{++}$-filter-regular element
for all $1\le i < t$ \cite {Tr2}.
\item[$\mathrm{(ii)}$] Trung and Verma in 2007  gave a version \cite[Theorem 1.4]{TV} of Theorem 3.3  with the help of sequences which they called $(\varepsilon_1,\ldots,\varepsilon_m)$-superficial sequences. Set
$$ \frak T  =\bigoplus_{n, n_1,\ldots,n_s\ge 0}\dfrac
{J^nI_1^{n_1}\cdots I_s^{n_s}}{J^{n+1}I_1^{n_1+1}\cdots
I_s^{n_s+1}}.$$  Recall  that
 $x$ is an {\it $i$-superficial element for $I_1,\ldots,I_s$} if $x \in I_i$ and the image $x^*$
 of $x$ in $\dfrac{I_i}{JI_1 \cdots I_{i-1}I_i^2I_{i+1}\cdots I_s}$
is an ${\frak T}_{++}$-filter-regular element in $\frak T$, i.e.,
$$(J^{n+1}I_1^{n_1+1}\cdots I_i^{n_i+2}\cdots I_s^{n_s+1}:x) \bigcap J^nI_1^{n_1}\cdots I_s^{n_s} = J^{n+1}I_1^{n_1+1}\cdots I_i^{n_i+1}\cdots I_s^{n_s+1}$$
for all $n, n_1,\ldots,n_s\gg0.$ And  if  $\varepsilon_1,\ldots,\varepsilon_m$ is a non-decreasing sequence of indices with $1\le \varepsilon_i\le s,$ then a sequence $x_1,\ldots,x_m$ is an  {\it $(\varepsilon_1,\ldots,\varepsilon_m)$-superficial sequence for $ I_1,\ldots,I_s$ } if for all $i= 1,\ldots,m$, $\bar{x}_i$ is an $\varepsilon_i$-superficial element for $\bar{J}, \bar{I}_1,\ldots,\bar{I}_s$, where $\bar{x}_i,\bar{J}, \bar{I}_1,\ldots,\bar{I}_s$ are the images of $x_i,J, I_1,\ldots,I_s$ in
$\frac{A}{(x_1,\ldots,x_{i-1})},$ respectively \cite{TV}.
 \end{itemize}
 \noindent
{\bf Remark 3.5.} In the case of two ideals $J, I_1,$  Manh in
\cite[Note 3.3.16(i)]{M} proved that sequences in Note 3.4(i) are
the sequences in Note 3.4(ii) when $s = 1$. Indeed,  it can be
verified that if $x \in I_1$ satisfies the conditions in Note
3.4(i), then
$$ (J^{m+1}I_1^{n+2}: x)\bigcap J^mI_1^{n+1} = J^{m+1}I_1^{n+1}\;
\mathrm{ and}\; (J^{m}I_1^{n+2}: x)\bigcap J^{m}I_1^{n} =
J^{m}I_1^{n+1} $$ for all $m, n \gg 0.$ Then for all $m, n \gg 0,$
we have
\begin{align*}
 (J^{m+1}I_1^{n+2}: x)\bigcap J^mI_1^{n}
 &=(J^{m+1}I_1^{n+2}: x)\bigcap (J^{m}I_1^{n+2}: x)\bigcap J^mI_1^{n}\\
 &= (J^{m+1}I_1^{n+2}: x)\bigcap J^mI_1^{n+1}
 = J^{m+1}I_1^{n+1}.
 \end{align*}
 \enlargethispage{1 cm}
 Thus  $x$ is an superficial element for $J, I_1$ as in Note 3.4(ii).

 Basing on \cite[Remark 4.1]{VM}, the authors of \cite{VM} constructed the following sequences that yield
filter-regular sequences in  $\bigoplus_{n_1,\ldots,n_s\ge 0}\dfrac{I_1^{n_1}\cdots I_s^{n_s}{M}}{I_1^{n_1+1}I_2^{n_2}\cdots I_s^{n_s}{M}}.$

\vskip 0.2cm \noindent {\bf Definition 3.6} \cite{VM}. Let $M$ be
a finitely generated  $A$-module. Let $\frak I,I_1,\ldots,I_s$ be
ideals such that $I=I_1\cdots I_s$  is not contained in $
\sqrt{\mathrm{Ann}{M}}$. An element $x \in A$ is called a {\it
superficial element of $(I_1,\ldots, I_s)$  with respect to $\frak
I$ and $M$} if there exists $i \in \{1, \ldots, s\}$ such that $x
\in I_i$ and

\begin{itemize}

    \item[$\mathrm{(i)}$]
 $({\frak I}I_1^{n_1}\cdots I_i^{n_i+1}\cdots I_s^{n_s}M:x) \bigcap I_1^{n_1}\cdots I_s^{n_s}M
 = {\frak I}I_1^{n_1}\cdots I_i^{n_i}\cdots I_s^{n_s}M$ \\for all  $ n_1,\ldots,n_s\gg0.$

\item[$\mathrm{(ii)}$] $xM\bigcap {I_1}^{n_1} \cdots I_i^{n_i+1}\cdots I_s^{n_s}M
= x{I_1}^{n_1}\cdots I_i^{n_i}\cdots I_s^{n_s}M$
for all $n_1,\ldots,n_s\gg0.$
\end{itemize}
Let $x_1, \ldots, x_t$ be a sequence in $A$. For each $0 \le i < t$, set ${M}_i = \dfrac{M}{(x_1, \ldots, x_{i})M}$. Then
$x_1, \ldots, x_t$\; is called a {\it superficial sequence  of $(I_1,\ldots, I_s)$ with respect to $\frak I$ and $M$} if $x_{i + 1}$ is a superficial element  of $({I}_1,\ldots, {I}_s)$  with respect to ${\frak I}$ and $M_i$ for $i = 0, \ldots, t - 1$.

The following main theorem  shows  that these sequences and
weak-(FC)-sequences are the same. \vskip 0.2cm \noindent {\bf
Theorem 3.7.}\;{\it  An element $x \in I_i$ is a superficial
element of $(I_1,\ldots, I_s)$
 with respect to $\frak I = I_1$ and $M$ if and only if  $x$  is a weak-$(FC)$-element of $M$ with respect to $(I_1,\ldots, I_s).$}
\begin{proof}\; Note that the condition (ii) in Definition 3.2 and the condition (ii) in Definition 3.6 are the same:
$xM\bigcap {I_1}^{n_1} \cdots I_i^{n_i+1}\cdots I_s^{n_s}M
= x{I_1}^{n_1}\cdots I_i^{n_i}\cdots I_s^{n_s}M$
for all $n_1,\ldots,n_s\gg0.$ Moreover, we get, by this equation and since $\frak I = I_1$ that
\begin{align*}
&(\frak I^kI_1^{n_1}\cdots I_i^{n_i+1}\cdots I_s^{n_s}M:x) \bigcap I_1^{n_1}\cdots I_s^{n_s}M \\
&= \big((\frak I^kI_1^{n_1}\cdots I_i^{n_i+1}\cdots I_s^{n_s}M \bigcap xM):x\big) \bigcap I_1^{n_1}\cdots I_s^{n_s}M \\
&= \big(x\frak I^kI_1^{n_1}\cdots I_i^{n_i}\cdots I_s^{n_s}M:x\big) \bigcap I_1^{n_1}\cdots I_s^{n_s}M \\
&= \big(\frak I^kI_1^{n_1}\cdots I_i^{n_i}\cdots I_s^{n_s}M + (0_M : x)\big) \bigcap I_1^{n_1}\cdots I_s^{n_s}M\\
&= \frak I^kI_1^{n_1}\cdots I_i^{n_i}\cdots I_s^{n_s}M + (0_M : x) \bigcap I_1^{n_1}\cdots I_s^{n_s}M
\end{align*}
 for all $k\geq 1$ and $ n_1,\ldots,n_s\gg0.$ Consequently
$$
 (\frak I^kI_1^{n_1}\cdots I_i^{n_i+1}\cdots I_s^{n_s}M:x) \bigcap I_1^{n_1}\cdots I_s^{n_s}M$$
$$\;\;\;\;\;\;\;\;\;\;=\frak I^kI_1^{n_1}\cdots I_i^{n_i}\cdots I_s^{n_s}M + (0_M : x) \bigcap I_1^{n_1}\cdots I_s^{n_s}M \eqno(1)
$$
for all $k\geq 1$ and $ n_1,\ldots,n_s\gg0.$
Let $x \in I_i$ be a superficial element of $(I_1,\ldots, I_s)$  with respect to $\frak I = I_1$ and $M$ as in Definition 3.6. Then
$$(\frak II_1^{n_1}\cdots I_i^{n_i+1}\cdots I_s^{n_s}M:x) \bigcap I_1^{n_1}\cdots I_s^{n_s}M  = \frak II_1^{n_1}\cdots I_i^{n_i}\cdots I_s^{n_s}M \eqno{(2)}$$
for all $ n_1,\ldots,n_s\gg0.$
We show by induction on $k \geq 1$ that
$$(\frak I^kI_1^{n_1}\cdots I_i^{n_i+1}\cdots I_s^{n_s}M:x) \bigcap I_1^{n_1}\cdots I_s^{n_s}M  = \frak I^kI_1^{n_1}\cdots I_i^{n_i}\cdots I_s^{n_s}M \eqno (3)$$
for all $k \geq 1$ and $ n_1,\ldots,n_s\gg0.$ The case that $k = 1$, (3) is true by (2). Assume that
$$(\frak I^kI_1^{n_1}\cdots I_i^{n_i+1}\cdots I_s^{n_s}M:x) \bigcap I_1^{n_1}\cdots I_s^{n_s}M  = \frak I^kI_1^{n_1}\cdots I_i^{n_i}\cdots I_s^{n_s}M $$
for all $ n_1,\ldots,n_s\gg0.$ It holds that
\begin{align*}
&(\frak I^{k + 1}I_1^{n_1}\cdots I_i^{n_i+1}\cdots I_s^{n_s}M:x) \bigcap I_1^{n_1}\cdots I_s^{n_s}M  \\
&=(\frak I^{k + 1}I_1^{n_1}\cdots I_i^{n_i+1}\cdots I_s^{n_s}M:x) \bigcap (\frak I^kI_1^{n_1}\cdots I_i^{n_i+1}\cdots I_s^{n_s}M:x) \bigcap I_1^{n_1}\cdots I_s^{n_s}M  \\
&=(\frak I^{k + 1}I_1^{n_1}\cdots I_i^{n_i+1}\cdots I_s^{n_s}M:x) \bigcap \frak I^kI_1^{n_1}\cdots I_s^{n_s}M  \\
&=(\frak II_1^{n_1+k}\cdots I_i^{n_i+1}\cdots I_s^{n_s}M:x) \bigcap I_1^{n_1 + k}\cdots I_s^{n_s}M  \\
&= \frak II_1^{n_1+k}\cdots I_i^{n_i}\cdots I_s^{n_s}M\\
&= \frak I^{k + 1}I_1^{n_1}\cdots I_i^{n_i}\cdots I_s^{n_s}M
\end{align*}
for all  $ n_1,\ldots,n_s\gg0.$
 Hence the induction is complete and we get (3).
Combining (3) with (1) we obtain
$$\frak I^kI_1^{n_1}\cdots I_i^{n_i}\cdots I_s^{n_s}M + (0_M : x) \bigcap I_1^{n_1}\cdots I_s^{n_s}M = \frak I^kI_1^{n_1}\cdots I_i^{n_i}\cdots I_s^{n_s}M$$
for all $k\geq 1$ and $ n_1,\ldots,n_s\gg0.$ Therefore
$$(0_M : x) \bigcap I_1^{n_1}\cdots I_s^{n_s}M \subset \frak I^kI_1^{n_1}\cdots I_i^{n_i}\cdots I_s^{n_s}M$$
for all $k\geq 1$ and $ n_1,\ldots,n_s\gg0.$ Hence
$$(0_M : x) \bigcap I_1^{n_1}\cdots I_s^{n_s}M \subset \bigcap_{k\geq 1}\frak I^kI_1^{n_1}\cdots I_i^{n_i}\cdots I_s^{n_s}M = 0$$
 for $ n_1,\ldots,n_s\gg0.$ Thus, $x$ is an $I$-filter-regular element with respect to $M$ and hence $x$ is a weak-(FC)-element of $M$ with respect to $(I_1,\ldots, I_s).$
Now, suppose that $x \in I_i$ is a weak-(FC)-element of $M$ with respect to $(I_1,\ldots, I_s).$  Since $x$ is an $I$-filter-regular element with respect to $M$, it follows that
$$(0_M: x) \bigcap I_1^{n_1}\cdots I_s^{n_s}M = 0
\;\;\mathrm{for \;\; all}\;\;  n_1,\ldots,n_s\gg0.$$
Indeed,   by Artin-Rees Lemma,  there exist positive integers
$u_1,\ldots,u_s$ such that $$\begin{array}{ll} I_1^{n_1}\cdots
I_s^{n_s}M\bigcap (0_M:I^\infty)&=I_1^{n_1-u_1}\cdots
I_s^{n_s-u_s}[ I_1^{u_1}\cdots I_s^{u_s}M\bigcap
(0_M:I^\infty)]\\&\subseteq I_1^{n_1-u_1}\cdots I_s^{n_s-u_s}
(0_M:I^\infty)\end{array}$$ for all $n_1\ge u_1,\ldots,n_s\ge
u_s$.  Since $M$  is noetherian, there exists a positive integer $u$ such that
$0_M:I^\infty=0_M:I^u.$  Therefore
$$\begin{array}{ll}I_1^{n_1-u_1}\cdots I_s^{n_s-u_s}
(0_M:I^\infty)&=I_1^{n_1-u_1}\cdots I_s^{n_s-u_s}
(0_M:I^u)\\&=I_1^{n_1-u_1}\cdots I_s^{n_s-u_s} (0_M:I_1^u\cdots
I_s^u)=0\end{array}$$ for all $n_1,\ldots,n_s\gg0.$ Remember that
$x$ is an $I$-filter-regular element with respect to $M,$ $0_M:x
\subseteq 0_M: I^{\infty}$. Thus
$$(0_M: x) \bigcap I_1^{n_1}\cdots I_s^{n_s}M \subset I_1^{n_1}\cdots
I_s^{n_s}M\bigcap (0_M:I^\infty) = 0$$ for all  $n_1,\ldots,n_s\gg0$. Hence
$$(\frak II_1^{n_1}\cdots I_i^{n_i+1}\cdots I_s^{n_s}M:x) \bigcap I_1^{n_1}\cdots I_s^{n_s}M = \frak II_1^{n_1}\cdots I_i^{n_i}\cdots I_s^{n_s}M$$
for all $ n_1,\ldots,n_s\gg0$ by (1). Thus $x$ satisfies the condition (i) in Definition 3.6 and hence $x$ is a superficial element of $(I_1,\ldots, I_s)$  with respect to $\frak I = I_1$ and $M$.
$\blacksquare$
\end{proof}
From the above facts, we obtain some  following comments.

\vskip 0.2cm
  \noindent
{\bf Remark 3.8.}   On the one hand, \cite[Remark 4.3 and Remark 4.7]{VM} showed that the superficial sequences of $\frak m$-primary ideals in Definition 3.1 and $(\varepsilon_1,\ldots,\varepsilon_m)$-superficial sequences in Note 3.4(ii)
are  the superficial sequences in  Definition 3.6. On the other hand,  sequences  in  Note 3.4(i) are sequences in Note 3.4 (ii) when $s=1$ by Remark 3.5. Hence sequences of $\mathfrak{m}$-primary ideals in Definition 3.1 and sequences in Note 3.4 are sequences in  Definition 3.6. But sequences in  Definition 3.6 and weak-(FC)-sequences in Definition 3.2 are the same by  Theorem 3.7.
 Hence the superficial sequences in \cite{Te, Tr2, TV, VM} are weak-(FC)-sequences.
From this it follows that  Theorem 3.3 \cite[Theorem 3.4]{Vi1} in
2000 covers the results of Risler   and  Teissier \cite{Te} in
1973; Trung \cite[Theorem 3.4]{Tr2} in 2001(see \cite {M}); Trung
and  Verma \cite[Theorem 1.4]{TV} in 2007 (see \cite{DV, DV1,
VT}). Moreover, \cite[Theorem 4.5]{VM} and \cite[Theorem 3.4]{Vi1}
are equivalent by Theorem 3.7.

  \vskip 0.2cm
  \noindent
{\bf Remark 3.9.} Theorem 3.7 and \cite[Remark 4.1]{VM}   seem to
account well for the minimum of conditions of weak-(FC)-sequences
that is used as a tool in interpreting  mixed multiplicities of
ideals into the Hilbert-Samuel multiplicity. This explains why the
superficial sequences in \cite{Te, Tr2, TV, VM} are
weak-(FC)-sequences. Note that weak-(FC)-sequences are also an
useful tool to study the multiplicity and the Cohen-Macaulayness
of blow-up rings (see e.g. \cite{DV2, Vi4, Vi5, Vi6, VT1, VT4}).

At this point we would like to emphasize that: On
the one hand the authors of \cite{Tr2,TV}
 used Viet's  results and method in \cite{Vi1} for their works with
 replacing weak-(FC)-sequences by other sequences that they called "filter-regular sequences".
 On the other hand one seemed to omit  Viet's work. This  causes confusion in citations,
 for instance, J. Huh cited
\cite[Theorem 3.4]{Vi1} in 2000 as Trung-Verma's theorem in 2007
(see \cite [Theorem 5]{JH}).

\end{document}